\documentclass[a4paper,fleqn]{cas-sc}

\usepackage[authoryear,longnamesfirst]{natbib}

\def\tsc#1{\csdef{#1}{\textsc{\lowercase{#1}}\xspace}}
\tsc{WGM}
\tsc{QE}
\tsc{EP}
\tsc{PMS}
\tsc{BEC}
\tsc{DE}
\begin{document}
\let\WriteBookmarks\relax
\def\floatpagepagefraction{1}
\def\textpagefraction{.001}
\shorttitle{}
\shortauthors{}


\title [mode = title]{Controllability and observability of linear multi-agent systems over matrix-weighted signed networks}                      
\tnotemark[1]
\tnotetext[1]{This work was supported by the National Natural Science Foundation of China (Grant nos. 61873136 and 62033007), Taishan Scholars Climbing Program of Shandong Province of China and Taishan Scholars Project of Shandong Province of China (No. ts20190930)}
\author[1,2]{Lanhao Zhao}
\address[1]{Institute of Complexity Science,College of Automation, Qingdao University, Qingdao 266071, China}
\author[1,2]{Zhijian Ji}
\cormark[1]
\ead{jizhijian@pku.edu.cn}
\address[2]{Shandong Key Laboratory of Industrial Control Technology, Qingdao 266071, China}
\cortext[cor1]{Corresponding author}
\author[3]{Yungang Liu}
\address[3]{School of Control Science and Engineering, Shandong University, Jinan, 250061, China}
\author[1,2]{Chong Lin}

\begin{abstract}
In this paper, the controllability and observability of linear multi-agent systems over matrix-weighted signed networks are analyzed. Firstly, the definition of equitable partition of matrix-weighted signed multi-agent system is given, and the upper bound of controllable subspace and a necessary condition of controllability are obtained by combining the restriction conditions of the coefficient matrix and matrix weight for the case of fixed and switching topologies, respectively. The influence of different selection methods of coefficient matrices on the results is discussed. Secondly, for the case of heterogeneous systems, the upper bound of controllable subspace and the necessary condition of controllability are obtained when the dynamics of individuals in the same cell are the same. Thirdly, sufficient conditions for controllable and uncontrollable union graphs are obtained by taking advantage of the concept of switched systems and equitable partitions, respectively. Finally, a necessary condition of observability is obtained in terms of the dual system and the constraints of the coefficient matrix, and the relationship between the observability and the controllability of the matrix-weighted signed multi-agent systems is discussed.
\end{abstract}


\maketitle
\newtheorem{theorem}{Theorem}
\newdefinition{lemma}{Lemma}
\newdefinition{rmk}{Remark}
\newproof{pf}{Proof}
\newdefinition{exm}{Example}
\newdefinition{de}{Definition}
\newdefinition{co}{Corollary}
\newdefinition{pro}{Proposition}

\section{Introduction}

Recently, the cooperative control of multiple individuals, such as UAV formation control and social network modeling, has attracted increasing attention of researchers. The concept of multi-agent system (MAS) was introduced to solve the problem of cooperative control. The so-called multi-agent system can be regarded as an ordered combination of multiple systems. Each individual updates their states through neighbor information, resulting in collective behaviors such as clustering and separation. Especially, many works of multi-agent systems had emerged, such as consensus \cite{1,2,3,4}, controllability \cite{5,6,7,8,9,10,11}, observability \cite{12,13,14}, stabilisability \cite{15}, optimal control \cite{16,17,18,19,20,21} and model-free adaptive control \cite{22}.
 
Controllability and observability are fundamental topics in control theory. The concept of controllability provides helps for the design of control law, and observability describes an ability to use partial information to figure out the whole system states. Kalman \cite{23} first proposed the concept of controllability, and it was generalized to multi-agent systems by Tanner \cite{24}. The system is controllable if the leader's state can be adjusted through an external input signal to guide the follower to any desired state. Most of the existing results of controllability and observability are obtained through the structural characteristics of graphs, especially the partition of graphs. Analyzing the controllability from the perspective of graph theory, especially partition, can directly determine the controllability of the system and save the calculation time. For example, the effects of network structure on system controllability from the perspective of equitable or distance partition (EP or DP) were discussed in \cite{25,26,27}, the results show that the upper and lower bounds of controllable subspace can be directly determined according to the number of cells and path length in the network, which provides convenience for the determination of system controllability. While \cite{28,29} discussed the structural controllability of different networks from the perspective of structured systems. It is worth noting that these works deal with the simplest cases. To make the controllability results proposed in these studies more applicable to various scenarios, some related work has appeared in recent years, such as Guan et al.\cite{30} extended the concept of equitable partition to the case of the directed graph. Ji et al discussed the case of delay or switching topology \cite{31}, signed graph \cite{32,33,34,35} was also discussed in recent years.  

Most of the existing results on the controllability and observability of multi-agent systems are mainly based on scalar-weighted networks. The concept of matrix weighted network is introduced to describe more complex systems, for example, opinion dynamics \cite{36} and formation control \cite{37}. Compared with the scalar-weighted networks, the matrix-weighted networks can manifest the interaction between the different states of individuals and the states of corresponding neighbors more precisely. It is non-trivial to analyze its controllability although the matrix weight network seems to be a natural extension of the scalar weight network. For instance, additional complexity is introduced because the rank of weight matrix is not fixed when describing the dimension of controllable subspace. \cite{38} analyzed the controllability of matrix-weighted network by using the method of almost equitable partition and distance partition and gave the upper and lower bounds of controllable subspace. Notably, the weight matrix is assumed to be positive definite or positive semi-definite in their work. However, the weight matrix may also be negative definite and negative semi-definite in more cases. The consensus problem under the assumption that the weight matrix can be negative definite or negative semi-definite was discussed in \cite{39}, however, this assumption has not been taken into account in the controllability analysis from the perspective of graph theory. Therefore, it is meaningful to consider the controllability and observability of the matrix-weighted signed networks under the assumption.

In addition, In the above work, the network topology is usually fixed, and each individual is described by a first-order integrator. In engineering practice, the network topology may change and the individual dynamics in the network are also different. This brings a challenge to analyze the controllability of the system by using graph theory. Firstly, the concept of equitable partition is extended to directed signed multi-agent systems with switched topology in \cite{40}, however, the concept of equitable partition is not used to analyze the relationship between union graph controllability and system controllability. Secondly, \cite{41} analyzes the controllability of multi-agent systems whose individual dynamics is a general linear system and reveals the new characteristics of the system compared with the first-order integrator system. Therefore, for matrix-weighted networks, it is meaningful to analyze the characteristics of individual dynamics when it is a general linear system. Finally, the difference of individuals dynamics in the network is also a challenge to the controllability analysis, especially for graphical methods. Specifically, most of the existing results of system controllability from the perspective of graph theory analyze the influence of some special structures based on the same dynamics of each individual. However, if we use different colors to label network nodes with different dynamics, the distribution law of nodes with different colors on the graph is a generalized graphic feature, which needs to be further studied in the conventional graphic feature. Therefore, it is necessary to analyze how individual dynamics affect the results of controllability and observability of the matrix-weighted signed networks. 

In this paper, the controllability and observability of linear multi-agent systems over matrix-weighted signed networks are analyzed. The main contributions of this paper are as follows. 

(1) The definition of equitable partition of matrix-weighted signed multi-agent system is first given, and the upper bound of controllable subspace and a necessary graphic condition of controllability are obtained by combining the restriction conditions of the coefficient matrix and matrix weight for the case of fixed and switching topologies, respectively. In particular, the equivalence of the controllability results based on the concept of equitable partition for first-order integrator systems and general linear systems under appropriate input matrix selection(it called complete control input) conditions is discussed.

(2) For the case of heterogeneous systems, the upper bound of controllable subspace and the necessary condition of controllability are obtained when the dynamics of individuals in the same cell are the same. To our knowledge, this is the first time to use the concept of equitable partition (EP) to deal with the controllability of heterogeneous systems whether for scalar or matrix weighted networks. The meaning of 'equitable' is extended by combining the characteristics of graph theory with the distribution law of each individual dynamics.

(3) Sufficient conditions for controllable and uncontrollable union graphs are obtained by taking advantage of the concept of switched systems and equitable partitions, respectively. As far as we know, this is the first time that the sufficient condition of uncontrollable union graph is given by using the concept of EP. 

(4) A necessary condition of observability is obtained in terms of the dual system and the constraints of the coefficient matrix, and the relationship between the observability and the controllability of the matrix-weighted signed multi-agent systems is discussed. The observability results in this paper are more applicable compared with the existing graph conditions of observability obtained under the conditions of the first-order integrator, undirected and unweighted graph.

The rest of this work is organized as follows. The notations and some preliminaries used in this paper are provided in Section 2. The model of matrix-weighted signed networks under general linear dynamics is established in Section 3. Section 4, 5, and 6 present the main results of controllability and observability of matrix weighted signed networks in different cases,  respectively. Some numerical examples are shown in Section 7, and the conclusions are summarized in Section 8.

\section{Notations and Preliminaries}

In this paper, $R$ stands for the set of real number, $I_{m}$ and $0_{m}$ denote identity matrix and zero matrix with dimension $m$, respectively. For a symmetric matrix $H \in \mathbb{R}^{n \times n}\left(H^{T}=H\right)$, it is denoted by $H \succ 0$ ($H \succeq 0$ ) when $H$ is positive definite (positive semi-definite), and $H$ is called negative definite (negative semi-definite) when $-H$ is positive definite (positive semi-definite). $|H|=H$ when $H \succ 0$ or $H \succeq 0$ and $|H|=-H$ when $H \prec 0$ or $H \preceq 0$.  $im(P)$ represents the vector space generated by the columns of matrix $P$. Graph theory: $G=\{V, E, \mathcal{A}\}$ is called as matrix-weighted signed graph, in which $V=\left\{v_{1}, v_{2} \cdots v_{n}\right\}$ represents the vertex set.  $\mathcal{A}=\left[\mathcal{A}_{i j}\right] \in R^{nd \times nd}$ is the weighted adjacency matrix which $\mathcal{A}_{i j}$ is a symmetric matrix belonging to $\mathbb{R}^{d \times d}$.  $E=E_{+}+E_{-}$ represents the edge set, where  ${E}_{+}=\left\{(j, i) \mid \mathcal{A}_{i j}\succ 0\cup \mathcal{A}_{i j}\succeq 0\right\}$ and ${E}_{-}=\left\{(j, i) \mid \mathcal{A}_{i j}\prec 0\cup  \mathcal{A}_{i j}\preceq 0\right\}$ denote the sets of positive and negative edges, respectively. The set of neighbors of agent $v_{i}$ is represented by $N_{i}=N_{i+}+N_{i-}$. $N_{i+}=\left\{j \in {V} \mid(j, i) \in {E}_{+}\right\}$ and  $N_{i-}=\left\{j \in {V} \mid(j, i) \in {E}_{-}\right\}$ represent the set of positive and negative neighbors, respectively. $d_{i}=\sum_{j\in{N}_{i}}\left|\mathcal{A}_{i j}\right|$ denotes the degree of $i$ for matrix-weighted signed graph. Let $L=D-\mathcal{A}$ be the Laplacian matrix of $G$, where $D=\operatorname{diag}\left(d_{1}, \cdots, d_{n}\right)$. The entries of $L$ can be writen as follows
$$
l_{i j}=\left\{\begin{array}{ll}
	d_{i}, & i=j \\
	-\mathcal{A}_{i j}, & i \neq j
\end{array}\right . 
$$

Graph partition: For the vertex set $V$ of a graph, its subset $V_{i}$ is called a cell. It is called a trivial cell  if the cell contains only one vertex, otherwise it is a nontrivial cell. If any vertex in $V_{1}$ also belongs to $V_{2}$,
then $V_{1}$ is a sub-cell of $V_{2}$. We define $\pi=\left\{V_{1}, V_{2}, \cdots, V_{k}\right\}$. Then $\pi$ is a partition of graph when $V_{i} \cap V_{j}=\varnothing$ and $\bigcup_{i} V_{i}=V$ for $0<i, j<k$ and $i \neq j$. The characteristic matrix is
$$
P_{i j} \triangleq\left\{\begin{array}{ll}
	I_{d\times d} , & i \in V_{j} \\
	0_{d\times d} , & i \notin V_{j}
\end{array}\right.
$$
\begin{exm}
	The characteristic matrix of $\pi=\{\{1,2\},\{3,4,5\}\}$ is
	$$
	P(\pi)=\left[\begin{array}{ll}
		I_{d\times d} & 0_{d\times d} \\
		I_{d\times d} & 0_{d\times d} \\
		0_{d\times d} & I_{d\times d} \\
		0_{d\times d} & I_{d\times d} \\
		0_{d\times d} & I_{d\times d}
	\end{array}\right].
	$$
\end{exm}
\section{Matrix-Weighted Signed Networks with General Linear Dynamics}
Consider a multiagent system with $n$ agents which state is denoted by the $x_{i}(t)\in \mathbb{R}^{d}$. The leader and follower are distinguished according to whether the agent receives external input signals. We assume that the first $m$ $(m<n)$ individuals are named as leaders where $V_{L}=\left\{v_{1}, \cdots, v_{m}\right\}$ is the leader set, and the follower set is represented by $V_{F}=\left\{v_{m+1}, \cdots, v_{n}\right\} $.

For the matrix-weighted signed networks under general linear dynamics, all followers are governed by the following dynamics
$$
\dot{x}_{i}(t)=A x_{i}(t)+B u_{i}(t).
$$

For every leader, its dynamics is
$$
\dot{x}_{i}(t)=A x_{i}(t)+B u_{i}(t)+C y_{i}(t)
$$
where $A \in R^{d \times d}, B \in R^{d \times p}, C \in R^{d \times q}$, $x_{i} \in R^{d}$ represents the state of the $ith$ agent, and $u_{i} \in R^{p}$ reflects the influence that each individual receives from the others, $y_{l} \in R^{q}$ represents an external input signal. The update rules based on neighbors are as follows
$$
u_{i}(k)=\sum_{j \in N_{i}} K\left[\mathcal{A}_{i j} x_{j}(k)-\left|\mathcal{A}_{i j}\right| x_{i}(k)\right]
$$
where $K$ stands for the feedback gain, $\mathcal{A}_{i j}$ is the connection weight between the individual $i$ and $j$. Denote $x(k)=\left[x_{1}(t), \cdots, x_{n}(t)\right]^{T}$ as the aggregate state vector and $y(k)=\left[y_{1}(t), \cdots, y_{l}(t)\right]^{T}$ as the control input vector. Then the expression of multiagent system can be written as
$$
\dot{x}(t)=\left(\left[\begin{array}{cccc}
	A & 0 & \cdots & 0 \\
	0 & A & \cdots & 0 \\
	\vdots & \vdots & \ddots & \vdots \\
	0 & 0 & \cdots & A
\end{array}\right]-\left[\begin{array}{cccc}
	B K & 0 & \cdots & 0 \\
	0 & B K & \cdots & 0 \\
	\vdots & \vdots & \ddots & \vdots \\
	0 & 0 & \cdots & B K
\end{array}\right] L\right) x(t)-\left(M\left[\begin{array}{cccc}
	C & 0 & \cdots & 0 \\
	0 & C & \cdots & 0 \\
	\vdots & \vdots & \ddots & \vdots \\
	0 & 0 & \cdots & C
\end{array}\right]\right) y(t)
$$
where matrix $M$ is used to distinguish the leaders from the followers.
$$
M_{i l}=\left\{\begin{array}{ll}
	I_{d\times d} & i=l \\
	0_{d\times d} & \text { otherwise }
\end{array}\right.
$$

For convenience, it can be rewritten as the following form
$$
\dot{x}(t)=\tilde{L} x(t)+\tilde{M} y(t)\label{1} \quad\quad\quad\quad (1)
$$
where
 $$
 \tilde{L}=\left(\left[\begin{array}{cccc}
 	A & 0 & \cdots & 0 \\
 	0 & A & \cdots & 0 \\
 	\vdots & \vdots & \ddots & \vdots \\
 	0 & 0 & \cdots & A
 \end{array}\right]-\left[\begin{array}{cccc}
 	B K & 0 & \cdots & 0 \\
 	0 & B K & \cdots & 0 \\
 	\vdots & \vdots & \ddots & \vdots \\
 	0 & 0 & \cdots & B K
 \end{array}\right] L\right), \tilde{M}=\left(M\left[\begin{array}{cccc}
 	C & 0 & \cdots & 0 \\
 	0 & C & \cdots & 0 \\
 	\vdots & \vdots & \ddots & \vdots \\
 	0 & 0 & \cdots & C
 \end{array}\right]\right).
 $$
 \begin{rmk}
 	System (1) includes both the influence of weight matrix and individual dynamics. Compared with the linear dynamics in scalar networks \cite{41}, the Kronecker product in the previous model of scalar networks is replaced by matrix product in the system (1). This makes some existing results based on Kronecker product properties may no longer be applicable especially when multiplication is needed. So this situation deserves further discussion for the new features brought by this matrix-weighted network.
 \end{rmk}
 
\section{Controllability and Graph partition}

In this section, the controllability of matrix-weighted signed networks is discussed by using graph partition. Firstly, the definition of equitable partition of matrix-weighted signed networks is given. Secondly, the controllability of matrix-weighted signed networks in fixed topology, switching topology and the heterogeneous system is analyzed by using the concept of equitable partition.

\begin{de}
    Denote a matrix-weighted signed graph as $G$, and let $\pi=\left\{V_{1}, V_{2}, \ldots, V_{k}\right\}$ be a partition of $G$. The partition $\pi$ is said to be an equitable partition (EP) of $G$ if for any $r, s \in V_{i}, i, j=1,2, \ldots, k$
	$$\sum_{t_{1} \in V_{j}, t_{1} \in {N}_{r+}} \mathcal{A}_{r t_{1}}=\sum_{t_{2} \in V_{j}, t_{2} \in {N}_{s+}} \mathcal{A}_{s t_{2}}$$
	and for any $r, s \in V_{i}, i, j=1,2, \ldots, k$
	$$\sum_{t_{1} \in  V_{j}, t_{1} \in N_{r-}} \mathcal{A}_{r t_{1}}=\sum_{t_{2} \in V_{j}, t_{2} \in \mathcal{N}_{s-}} \mathcal{A}_{s t_{2}}$$
	where $\left(t_{1}, r\right),\left(t_{2}, s\right) \in E$. 	
\end{de}
\begin{rmk}
	Compared with the definition for scaler weighted network in \cite{34}, The sum of scalar weights in previous definition is replaced by the sum of matrices. This shows the difference between matrix weighted network and scalar weighted network when using the concept of equitable partition to analyze controllability. 
\end{rmk}
\begin{rmk}
	Compared with the case that all the matrix weights are positive definite in \cite{38}, Definition 1 applies in the case that the matrix weights are negative definite or semi negative definite. Therefore, this extension is necessary and meaningful.
\end{rmk}

Denote $d(v_{i},Q)=\sum_{v_{j}\in Q}\left|\mathcal{A}_{i j}\right|$ and $d(V_{i},Q)=d(v,Q)$ for all $v \in V_{i}$, then we give the concept of quotient graph
\begin{de}
	For an equitable partition $\pi=$ $\left\{V_{1}, V_{2}, \ldots, V_{s}\right\}$ of a matrix-weighted network $G$, the quotient graph of $G$ over $\pi$ is a matrix-weighted network denoted by $G/\pi$ with the node set $V(G/\pi)=\left\{v_{1}, v_{2}, \ldots, v_{s}\right\}$ and the edge set is $E(G/ \pi)=\left\{\left(v_{i}, v_{j}\right) \mid d\left(V_{i}, V_{j}\right) \neq 0_{d \times d}\right\},$ where
	the weight of edge $\left(V_{i}, V_{j}\right)$ is $d\left(V_{i}, V_{j}\right)$ for $i \neq j \in \underline{s}$.
\end{de}

Denote $L_{\pi}$ as the Laplacian matrix of $G/\pi$
$$
\left(L_{\pi}\right)_{i j}=\left\{\begin{array}{ll}
	\sum_{v_{j} \in V(G/\pi)} d\left(V_{i}, V_{j}\right), & i=j \\
	-d\left(V_{i}, V_{j}\right), & i \neq j
\end{array}\right.
$$
\begin{lemma}
	$ \pi=\left\{V_{1}, V_{2}, \cdots, V_{s}\right\}$ is an EP for matrix-weighted signed graph $G$ and $P_{\pi}$ is the characteristic matrix. Then $L$ satisfies
	$$
	L P_{\pi}=P_{\pi} L_{\pi}.
	$$
	Furthermore $i m\left(P_{\pi}\right)$ is $L-$ invariant. 	
\end{lemma}
\begin{pf}
	The proof is similar to Lemma 1 in \cite{38}, and thus is omitted.
\end{pf}

\begin{rmk}
	 Lemma 1 shows the relationship between the characteristic matrix and the Laplacian matrix. For the first-order system, the upper bound of the controllable subspace can be obtained by the expression of the controllable subspace in the geometric theory of linear systems \cite{38}. However, how to deal with the unknown coefficient matrix is a challenge for general linear dynamics and worth further discussion.
\end{rmk}

\subsection{The case of fixed topology}
We first discuss the case of fixed topologies. The so-called fixed topology means that the network topology between individuals does not change in the evolution process of multi-agent systems. For system (1), the system is controllable if and only if the matrix $\left[\tilde{M} \quad\tilde{L} \tilde{M} \quad\tilde{L}^{2} \tilde{M} \cdots \tilde{L}^{dn-1} \tilde{M}\right]$ has full row rank, and the controllable subspace of system (1) is
$$
\mathcal{W}=\left\langle \tilde{L} \mid \tilde{M}\right\rangle=i m(\tilde{M})+\tilde{L} \times i m(\tilde{M})+\cdots+\tilde{L}^{dn-1} \times i m(\tilde{M}).
$$
It is a minimal $\tilde{L}-$ invariant subspace containing $\operatorname{im}(\tilde{M})$.

Denote
$$
\tilde{P}_{\pi}=P_{\pi} \cdot\left[\begin{array}{cccc}
	C & 0 & \cdots & 0 \\
	0 & C & \cdots & 0 \\
	\vdots & \vdots & \ddots & \vdots\\
	0 & 0 & \cdots & C
\end{array}\right].
$$
 
\begin{lemma}
	$ \pi=\left\{V_{1}, V_{2}, \cdots, V_{s}\right\}$ is an EP for matrix-weighted graph $G$ and $P_{\pi}$ is the characteristic matrix. 
If the matrix 
	$$
	Q=\left[\begin{array}{cccc}
		Q_{1} & 0 & \cdots & 0 \\
		0 & Q_{1} & \cdots & 0 \\
		\vdots & \vdots & \ddots & \vdots \\
		0 & 0 & \cdots & Q_{1}
	\end{array}\right]-\left[\begin{array}{cccc}
		Q_{11} & Q_{12} & \cdots & Q_{1 \mathrm{~s}} \\
		Q_{21} & Q_{12} & \cdots & Q_{2 s} \\
		\vdots & \vdots & \ddots & \vdots \\
		Q_{n 1} & Q_{n 2} & \cdots & Q_{n s}
	\end{array}\right]
	$$
	exists, then there are matrices $Q_{1}, Q_{ij}$ that make  $A C = C Q_{1}$ and $B K L_{\pi ij} C = C Q_{ij}, (0 \leq i \leq n, 0 \leq j \leq s) $. It follows that
	$$
	\tilde{L} \tilde{P}_{\pi}=\tilde{P}_{\pi} Q.
	$$ 
	Furthermore $i m\left(\tilde{P}_{\pi}\right)$ is $\tilde{L}-$ invariant. 	
\end{lemma}
\begin{pf}
	$$
	\tilde{L}=\left(\left[\begin{array}{cccc}
		A & 0 & \cdots & 0 \\
		0 & A & \cdots & 0 \\
		\vdots & \vdots & \ddots & \vdots \\
		0 & 0 & \cdots & A
	\end{array}\right]-\left[\begin{array}{cccc}
		B K & 0 & \cdots & 0 \\
		0 & B K & \cdots & 0 \\
		\vdots & \vdots & \ddots & \vdots\\
		0 & 0 & \cdots & B K
	\end{array}\right] L\right), \tilde{P}_{\pi}=P_{\pi} \cdot\left[\begin{array}{cccc}
	C & 0 & \cdots & 0 \\
	0 & C & \cdots & 0 \\
	\vdots & \vdots & \ddots & \vdots \\
	0 & 0 & \cdots & C
\end{array}\right].
	$$
	According to Lemma 1
	$$
	\begin{aligned}
		\tilde{L} \tilde{P}_{\pi}&=\left[\begin{array}{cccc}
			A & 0 & \cdots & 0 \\
			0 & A & \cdots & 0 \\
			\vdots & \vdots & \ddots & \vdots \\
			0 & 0 & \cdots & A
		\end{array}\right] \cdot P_{\pi} \cdot\left[\begin{array}{cccc}
			C & 0 & \cdots & 0 \\
			0 & C & \cdots & 0 \\
			\vdots & \vdots & \ddots & \vdots \\
			0 & 0 & \cdots & C
		\end{array}\right]-\left[\begin{array}{cccc}
			B K & 0 & \cdots & 0 \\
			0 & B K & \cdots & 0 \\
			\vdots & \vdots & \ddots & \vdots\\
			0 & 0 & \cdots & B K
		\end{array}\right] \cdot L \cdot P_{\pi} \cdot\left[\begin{array}{cccc}
			C & 0 & \cdots & 0 \\
			0 & C & \cdots & 0 \\
			\vdots & \vdots & \ddots & \vdots \\
			0 & 0 & \cdots & C
		\end{array}\right]\\
		&=\left[\begin{array}{cccc}
			A P_{11} C & A P_{12} C & \cdots & A P_{1 s} C \\
			A P_{21} C & A P_{22} C & \cdots & A P_{2 s} C \\
			\vdots & \vdots & \ddots & \vdots \\
			A P_{n 1} C & A P_{n 2} C & \cdots & A P_{n s} C
		\end{array}\right]-\left[\begin{array}{cccc}
			B K & 0 & \cdots & 0 \\
			0 & B K & \cdots & 0 \\
			\vdots & \vdots & \ddots & \vdots \\
			0 & 0 & \cdots & B K
		\end{array}\right] \cdot P_{\pi} \cdot L_{\pi} \cdot\left[\begin{array}{cccc}
			C & 0 & \cdots & 0 \\
			0 & C & \cdots & 0 \\
			\vdots & \vdots & \ddots & \vdots\\
			0 & 0 & \cdots & C
		\end{array}\right]\\
		&=\left[\begin{array}{cccc}
			A P_{11} C & A P_{12} C & \cdots & A P_{1 s} C \\
			A P_{21} C & A P_{22} C & \cdots & A P_{2 s} C \\
			\vdots & \vdots & \ddots & \vdots \\
			A P_{n 1} C & A P_{n 2} C & \cdots & A P_{n s} C
		\end{array}\right]-\left[\begin{array}{cccc}
			B K P_{11} & B K P_{12} & \cdots & B K P_{1 s} \\
			B K P_{21} & B K P_{22} & \cdots & B K P_{2 s} \\
			\vdots & \vdots & \ddots & \vdots \\
			B K P_{n 1} & B K P_{n 2} & \cdots & B K P_{n s}
		\end{array}\right] \cdot L_{\pi} \cdot\left[\begin{array}{cccc}
			C & 0 & \cdots & 0 \\
			0 & C & \cdots & 0 \\
			\vdots & \vdots & \ddots & \vdots \\
			0 & 0 & \cdots & C
		\end{array}\right]
	\end{aligned}
	$$
	The specific expression of the characteristic matrix is
	$$
	P_{i j} \triangleq\left\{\begin{array}{ll}
		I_{d\times d} & i \in V_{j} \\
		0_{d\times d} & i \notin V_{j}
	\end{array}\right.
	$$
	
	Therefore, if $A C = C Q_{1}$ holds, then
	$$
	\left[\begin{array}{cccc}
		A P_{11} C & A P_{12} C & \cdots & A P_{1 s} C \\
		A P_{21} C & A P_{22} C & \cdots & A P_{2 s} C \\
		\vdots & \vdots & \ddots & \vdots \\
		A P_{n 1} C & A P_{n 2} C & \cdots & A P_{n s} C
	\end{array}\right]=P_{\pi} \cdot\left[\begin{array}{cccc}
		C & 0 & \cdots & 0 \\
		0 & C & \cdots & 0 \\
		\vdots & \vdots & \ddots & \vdots \\
		0 & 0 & \cdots & C
	\end{array}\right] \cdot\left[\begin{array}{cccc}
		Q_{1} & 0 & \cdots & 0 \\
		0 & Q_{1} & \cdots & 0 \\
		\vdots & \vdots & \ddots & \vdots \\
		0 & 0 & \cdots & Q_{1}
	\end{array}\right]
	$$
	
	And if $B K L_{\pi ij} C = C Q_{ij}, (0 \leq i \leq n, 0 \leq j \leq s) $ holds, then
	$$
	\left[\begin{array}{cccc}
		B K P_{11} & B K P_{12} & \cdots & B K P_{1 s} \\
		B K P_{21} & B K P_{22} & \cdots & B K P_{2 s} \\
		\vdots & \vdots & \ddots & \vdots \\
		B K P_{n 1} & B K P_{n 2} & \cdots & B K P_{n s}
	\end{array}\right] \cdot L_{\pi} \cdot\left[\begin{array}{cccc}
		C & 0 & \cdots & 0 \\
		0 & C & \cdots & 0 \\
		\vdots & \vdots & \ddots & \vdots\\
		0 & 0 & \cdots & C
	\end{array}\right]=P_{\pi} \cdot\left[\begin{array}{cccc}
		C & 0 & \cdots & 0 \\
		0 & C & \cdots & 0 \\
		\vdots & \vdots & \ddots & \vdots \\
		0 & 0 & \cdots & C
	\end{array}\right] \cdot\left[\begin{array}{cccc}
		Q_{11} & Q_{12} & \cdots & Q_{1 s} \\
		Q_{21} & Q_{22} & \cdots & Q_{2 s} \\
		\vdots & \vdots & \ddots & \vdots \\
		Q_{n 1} & Q_{n 2} & \cdots & Q_{n s}
	\end{array}\right]
	$$
	
	Denote $Q$ as
	$$
	Q=\left[\begin{array}{cccc}
		Q_{1} & 0 & \cdots & 0 \\
		0 & Q_{1} & \cdots & 0 \\
		\vdots & \vdots & \ddots & \vdots\\
		0 & 0 & \cdots & Q_{1}
	\end{array}\right]-\left[\begin{array}{cccc}
		Q_{11} & Q_{12} & \cdots & Q_{1 \mathrm{~s}} \\
		Q_{21} & Q_{12} & \cdots & Q_{2 s} \\
		\vdots & \vdots & \ddots & \vdots \\
		Q_{n 1} & Q_{n 2} & \cdots & Q_{n s}
	\end{array}\right]
	$$
	
	It can be derived that 
	$$
	\tilde{L} \tilde{P}_{\pi}=\tilde{P}_{\pi} Q
	$$
	if $A C = C Q_{1}$ and $B K L_{\pi ij} C = C Q_{ij}, (0 \leq i \leq n, 0 \leq j \leq s) $. Furthermore $i m\left(\tilde{P}_{\pi}\right)$ is $\tilde{L}-$ invariant.
\end{pf}

\begin{theorem}
	The controllable subspace $\mathcal{W}$ of the system (1) satisfies $\mathcal{W} \subseteq i m\left(\tilde{P}_{\pi}\right)$ if the matrix $Q$ exists.
\end{theorem}
\begin{pf}
	Every column of $M$ is also a column of $P_{\pi_{L}^{*}},$ and then every column of $\tilde{M}$ is also a column of $\tilde{P}_{\pi}$. It follows that $ im(\tilde{M}) \subseteq im\left(\tilde{P}_{\pi_{L}^{*}}\right)$. By Lemma $2$,  $\operatorname{im}(\tilde{P})$ is
	$\tilde{L}-$ invariant. Then
	$$
	\begin{aligned}
		\mathcal{W} &=im(\tilde{M})+\tilde{L} \times i m(\tilde{M})+\cdots+\tilde{L}^{dn-1} \times i m(\tilde{M}) \\
		& \subseteq im\left(\tilde{P}_{\pi}\right)+\tilde{L} \times i m\left(\tilde{P}_{\pi}\right)+\cdots+\tilde{L}^{dn-1} \times im\left(\tilde{P}_{\pi}\right) \\
		&=im\left(\tilde{P}_{\pi}\right).
	\end{aligned}
	$$
\end{pf}
\begin{co}
	Suppose that the matrix $Q$ exists, then every cell in the EP is trivial when system (1) is controllable.  
\end{co}
\begin{pf}
	If a certain cell is a nontrivial cell under controllable system (1), the dimension of $im\left(\tilde{P}_{\pi}\right)$ is less than
    $d n$. By Theorem $1,$ the dimension of the controllable subspace is also less than $dn,$ and accordingly system  (1) is uncontrollable. This is contrary to the assumption. Therefore, when system (1) is controllable, each cell is a trivial cell. The proof is complete. 
\end{pf}
\begin{rmk}
   Theorem 1 and Corollary 1 are similar to the results under the scalar weight network. For instance, $\tilde{L} \tilde{P}_{\pi}=\tilde{P}_{\pi} Q$ holds for matrix weighted networks base on Lemma 2. It is obtained that  $L P=P L_{\pi}$ if matrices $A$ is chosen as zero matrices and $B$, $K$ and $C$ are all chosen as identity matrices, which is same as $L P=P L_{\pi}$ for scaler weighted network \cite{34}. It also shows that the scalar weighted network is a special case of matrix weight network.
\end{rmk}
\begin{rmk}
	It is worth noting that Theorem 1 and Corollary 1 are also different from the conclusion under the scalar weighted network. It mainly in the following two points. The weight matrix can not be ignored and does not meet the commutative law in matrix operation compared with scalar weight. And the dimension of the weight matrix also needs to be considered when estimating the upper bound of the controllable subspace of the system. 
\end{rmk}

From Theorem 1 and Corollary 1, matrices $Q_{1}$ and $Q_{ij}$ need to be found to make $A C = C Q_{1}$ and $B K L_{\pi ij} C = C Q_{ij}, (0 \leq i \leq n, 0 \leq j \leq s) $ hold. There is no such limitation for the similar results in \cite{38}, so it is necessary to further discuss this condition.

First, we discuss the simplest case that matrices $A$, $B$, $K$, $L_{\pi i j}$ and $C$ are all chosen as identity matrices, and matrices $Q_{1}$ and $Q_{ij}$ can be found to make $A C = C Q_{1}$ and $B K L_{\pi ij} C = C Q_{ij}, (0 \leq i \leq n, 0 \leq j \leq s) $ hold for this case.
$$
Q_{1}=Q_{ij}=I_{d\times d}, (0 \leq i \leq n, 0 \leq j \leq s). 
$$
Then the conditions of Theorem 1 and Corollary 1 both hold in this case.

Furthermore, the matrix $L_{\pi ij}$ is not always chosen as the identity matrix in practice. So we discuss a more complicated case that matrices $A$, $B$, $K$ and $C$ are all chosen as identity matrices. Note that $Q_{1}$ and $Q_{ij}$ can still be found to make $A C = C Q_{1}$ and $B K L_{\pi ij} C = C Q_{ij}, (0 \leq i \leq n, 0 \leq j \leq s) $ hold for this case.
$$
Q_{1}=I_{d\times d}, Q_{ij}=L_{\pi ij}, (0 \leq i \leq n, 0 \leq j \leq s). 
$$
Then the conditions of Theorem 1 and Corollary 1 both hold in this case. It should be noted that the model considered in this case is the same to the one in \cite{38}.

In practical engineering, especially when individual dynamics is more complex, matrices $A$, $B$, and $K$ are rarely selected as identity matrices, and only matrix $C$ can be designed by our goals. At this point, if matrix $C$ is chosen as an identity matrix, matrices $Q_{1}$ and $Q_{ij}$ can still be found to make $A C = C Q_{1}$ and $B K L_{\pi ij} C = C Q_{ij}, (0 \leq i \leq n, 0 \leq j \leq s)$ hold for this case.
$$
Q_{1}=A, Q_{ij}=B K L_{\pi ij}, (0 \leq i \leq n, 0 \leq j \leq s). 
$$
The conditions of Theorem 1 and Corollary 1 both hold in this case.

Furthermore, if matrix $C$ is choosen as an invertible matrix, the matrices $Q_{1}$ and $Q_{ij}$ can still be found to make $A C = C Q_{1}$ and $B K L_{\pi ij} C = C Q_{ij}, (0 \leq i \leq n, 0 \leq j \leq s) $ hold for this case.
$$
Q_{1}=C^{-1} A C, Q_{ij}=C^{-1} B K L_{\pi ij} C, (0 \leq i \leq n, 0 \leq j \leq s). 
$$
Then the conditions of Theorem 1 and Corollary 1 both hold in this case.

Based on the above discussion, it can be obtained that matrix $Q$ is always exists when matrix $C$ is chosen as an invertible matrix. Therefore, this situation is worth discussing. We take C as the identity matrix for the convenience. When matrix C is the identity matrix, it means that the components of the external control signal correspond to the leader's state one by one. Then the results of general linear system and first-order integrator system are the same for the controllability analysis based on the concept of equitable partition, which is not affected by system dynamics. This also explains why there is no such limitation for the similar results in \cite{38}, and Theorem 1 is an extension of previous conclusions in \cite{38}.

Inspired by this situation, the concept of complete control input is introduced

\begin{de}
	For leaders which are described as
	$$
	\dot{x}_{i}(t)=A x_{i}(t)+B u_{i}(t)+C y_{i}(t),
	$$
	if $C$ is an invertible matrix, then the control input  $y_{i}(t)$ is called the complete control input.
\end{de}

According to Theorem 1, Corollary 1 and the Definition 3, the following result can be obtained
\begin{co}
	Every cell in the EP is trivial when system (1) is controllable with a complete control input.
\end{co}
\begin{rmk}
	In practical systems, matrix $C$ is sometimes not designable and not invertible. In this case, matrices $Q_{1}$ and $Q_{ij}$ need to be found to make $A C = C Q_{1}$ and $B K L_{\pi ij} C = C Q_{ij}, (0 \leq i \leq n, 0 \leq j \leq s) $ hold. it reflects the influence of the selection of coefficient matrix on controllability analysis.
\end{rmk}

\subsection{The case of switching topology}
In practical application, due to the influence of formation change, fault, obstacle avoidance, and the other factors, the connection relationship between individuals in multi-agent systems maybe change. The concept of switching topologies is introduced to deal with this situation. For the case of switching topology, denote $\left\{G_{1}, G_{2}, \ldots, G_{t}\right\}$ as the set of  switching topology. Function $\sigma(x): R^{+} \rightarrow\{1,2, \ldots, t\}$ is the switching signal which is used to determine the topological structure of the system at a certain time point. For convenience, we only consider the situation of fixed leaders during the switching process. Then the system can be written as 
$$
\dot{x}(t)=\tilde{L}_{\sigma(x)} x(t)+\tilde{M} y(t)\label{1} \quad\quad\quad\quad (2)
$$
where
$$
\tilde{L}_{\sigma(x)}=\left(\left[\begin{array}{cccc}
	A & 0 & \cdots & 0 \\
	0 & A & \cdots & 0 \\
	\vdots & \vdots & \ddots & \vdots\\
	0 & 0 & \cdots & A
\end{array}\right]-\left[\begin{array}{cccc}
	B K & 0 & \cdots & 0 \\
	0 & B K & \cdots & 0 \\
	\vdots & \vdots & \ddots & \vdots \\
	0 & 0 & \cdots & B K
\end{array}\right] L_{\sigma(x)}\right), \tilde{M}=\left(M\left[\begin{array}{cccc}
	C & 0 & \cdots & 0 \\
	0 & C & \cdots & 0 \\
	\vdots & \vdots & \ddots & \vdots \\
	0 & 0 & \cdots & C
\end{array}\right]\right).
$$

The switching sequence of the system (2) is
$$
\pi=\left\{\left(i_{1}, h_{1}\right),\left(i_{2}, h_{2}\right), \ldots,\left(i_{M}, h_{M}\right)\right\}
$$
where $i_{m} \in\{1,2, \ldots, t\}$, and $h_{m}=t_{m}-t_{m-1}>$ 0 represents the duration of the $m$-th topology. If for any initial state, there are one switch sequence and external inputs that take the states to zeros, then the system (2) is said to be controllable. The controllable state set is denoted as $\mathcal{C}(\pi)$ and the subspace sequence \cite{42} is 
$$
\mathbb{W}_{1}=\sum_{\sigma(x)=1}^{t}\left\langle \tilde{L}_{\sigma(x)} \mid \tilde{M}\right\rangle, \quad \mathbb{W}_{i+1}=\sum_{\sigma(x)=1}^{t}\left\langle\tilde{L}_{\sigma(x)} \mid \mathbb{W}_{i}\right\rangle, \quad i=1,2, \ldots
$$
we denote $\pi_{i}=\left\{V_{1}, V_{2}, \cdots, V_{s}\right\}$ as an EP for matrix-weighted signed graph $G_{i}$ and $\pi_{min}=\left\{V_{1}, V_{2}, \cdots, V_{s}\right\}$ is the least upper bound for all $\pi_{i}$. 

It can be seen that the expression of subspace is similar to that of fixed topology, so we can get the following results by using the concept of EP.
\begin{theorem}
For system (2), if the matrix $Q$ always exists for all $\tilde{L}_{i}$ corresponding to $G_{i}$, the dimension of its controllable state set $\mathcal{C}$ satisfies
$$
\operatorname{dim}(\mathcal{C}) \leq \operatorname{card}\left(\pi_{min}\right)\cdot d
$$
\end{theorem}
\begin{pf}
Since $\mathcal{C}=\mathbb{W}_{N}$ \cite{40}, the conclusion is equivalent to $\operatorname{dim}\left(\mathbb{W}_{N}\right) \leq$ $\operatorname{card}\left(\pi_{min}\right)\cdot d$. According to Lemma 2, if  $Q$ always exists for all $\tilde{L}_{i}$ corresponding to $G_{i}$, so for any $\pi_{i}$, $i \in\{1,2, \ldots, t\}, \operatorname{im}\left(P_{i}\right)$ is $\tilde{L}_{i}$ -invariant, where $P_{i}$ is the characteristic matrix of $\pi_{i}$.
$$
\begin{aligned}
	\mathbb{W}_{1} &=\sum_{i=1}^{t} \sum_{m=0}^{N-1}\left(\tilde{L}_{i}\right)^{m} \operatorname{im}\left(\tilde{M}\right) \\
	& \subseteq \sum_{i=1}^{t} \sum_{m=0}^{N-1}\left(\tilde{L}_{i}\right)^{m} \operatorname{im}\left(P_{i}\right)=\operatorname{im}\left(P_{i}\right) \\
	\mathbb{W}_{2} &=\sum_{i=1}^{t} \sum_{m=0}^{N-1}\left(\tilde{L}_{i}\right)^{m} \mathbb{W}_{1} \subseteq \sum_{i=1}^{t} \sum_{m=0}^{N-1}\left(\tilde{L}_{i}\right)^{m} \operatorname{im}\left(P_{i}\right) =\operatorname{im}\left(P_{i}\right)  \\
	& \ldots \\
	\mathbb{W}_{N} &=\sum_{i=1}^{t} \sum_{m=0}^{N-1}\left(\tilde{L}_{i}\right)^{m} \mathbb{W}_{N-1} \\
	& \subseteq \sum_{i=1}^{t} \sum_{m=0}^{N-1}\left(\tilde{L}_{i}\right)^{m} \operatorname{im}\left(P_{i}\right) =\operatorname{im}\left(P_{i}\right) 
\end{aligned}
$$
That is to say, $\operatorname{dim}\left(\mathbb{W}_{N}\right) \leq \operatorname{dim}\left(\operatorname{im}\left(P_{i}\right)\right)=\operatorname{card}(P_{i})\cdot d$, and then 
$$
\operatorname{dim}\left(\mathbb{W}_{N}\right) \leq \operatorname{card}\left(\pi_{min}\right)\cdot d.
$$

This completes the proof.
\end{pf}

\begin{rmk}
Theorem 2 shows that the dimension of the controllable subspace is limited by the nontrivial cell contained in each subgraph. This reflects the influence of switching topology on the controllability analysis based on the concept of equitable partition.
\end{rmk}

\subsection{The case of heterogeneous systems}
In practical applications, the models of different individuals are often in different due to the requirements of tasks or manufacturing reasons. We call this situation as heterogeneous system, which brings challenges to controllability analysis, especially from the perspective of the graphics method. Here, we consider matrix weighted signed networks, in which each individual can be described by general linear dynamics, and the coefficient matrix can be different. All followers take the following dynamics
$$
\dot{x}_{i}(t)=A_{i} x_{i}(t)+B_{i} u_{i}(t)
$$
and for every leader, its dynamics is
$$
\dot{x}_{i}(t)=A_{i} x_{i}(t)+B_{i} u_{i}(t)+C y_{i}(t)
$$
where $A_{i} \in R^{d \times d}, B_{i} \in R^{d \times p}, C \in R^{d \times q}$. Then the dynamics of multi-agent system can be written as
$$
\dot{x}(t)=\left(\left[\begin{array}{cccc}
	A_{1} & 0 & \cdots & 0 \\
	0 & A_{2} & \cdots & 0 \\
	\vdots & \vdots & \ddots & \vdots \\
	0 & 0 & \cdots & A_{n}
\end{array}\right]-\left[\begin{array}{cccc}
	B_{1} K L_{11} &B_{1} K L_{12} & \cdots & B_{1} K L_{1n} \\
	B_{2} K L_{21} & B_{2} K L_{22} & \cdots & B_{2} K L_{2n} \\
	\vdots & \vdots & \ddots & \vdots \\
	B_{n} K L_{n1} & B_{n} K L_{n2} & \cdots & B_{n} K L_{nn}
\end{array}\right]\right) x(t)-\left(M\left[\begin{array}{cccc}
	C & 0 & \cdots & 0 \\
	0 & C & \cdots & 0 \\
	\vdots & \vdots & \ddots & \vdots \\
	0 & 0 & \cdots & C
\end{array}\right]\right) y(t)
$$

For convenience, it can be rewritten as the following form
$$
\dot{x}(t)=\tilde{L}^{'} x(t)+\tilde{M} y(t)\label{1} \quad\quad\quad\quad (3)
$$
where
$$
\tilde{L}^{'}=\left(\left[\begin{array}{cccc}
	A_{1} & 0 & \cdots & 0 \\
	0 & A_{2} & \cdots & 0 \\
	\vdots & \vdots & \ddots & \vdots \\
	0 & 0 & \cdots & A_{n}
\end{array}\right]-\left[\begin{array}{cccc}
	B_{1} K L_{11} &B_{1} K L_{12} & \cdots & B_{1} K L_{1n} \\
	B_{2} K L_{21} & B_{2} K L_{22} & \cdots & B_{2} K L_{2n} \\
	\vdots & \vdots & \ddots & \vdots \\
	B_{n} K L_{n1} & B_{n} K L_{n2} & \cdots & B_{n} K L_{nn}
\end{array}\right]\right).
$$

To deal with different coefficient matrices, we assume that the dynamics of nodes in each nontrivial cell is the same, which means that $A_{i}=A_{j},B_{i}=B_{j}$ if nodes $i,j $ belong to the same cell, and the dynamics of nodes in different cells can be different. Then we can get one result similar to Lemma 2.

\begin{lemma}
	$ \pi=\left\{V_{1}, V_{2}, \cdots, V_{k}\right\}$ is an EP for matrix-weighted graph $G$ and $P_{\pi}$ is the corresponding characteristic matrix. Suppose $A_{i}=A_{j},B_{i}=B_{j}$ if nodes $i,j $ belong to same cell. When the matrix 
	$$
	Q^{'}=\left[\begin{array}{cccc}
		Q^{'}_{1} & 0 & \cdots & 0 \\
		0 & Q^{'}_{2} & \cdots & 0 \\
		\vdots & \vdots & \ddots & \vdots \\
		0 & 0 & \cdots & Q^{'}_{s}
	\end{array}\right]-\left[\begin{array}{cccc}
		Q^{'}_{11} & Q^{'}_{12} & \cdots & Q^{'}_{1 s} \\
		Q^{'}_{21} & Q^{'}_{12} & \cdots & Q^{'}_{2 s} \\
		\vdots & \vdots & \ddots & \vdots \\
		Q^{'}_{n 1} & Q^{'}_{n 2} & \cdots & Q^{'}_{n s}
	\end{array}\right]
	$$
	exists, which means that there are matrices $Q^{'}_{1}, Q^{'}_{ij}$ that make $A_{i} C = C Q^{'}_{i}$ and $B_{i} K L_{\pi ij} C = C Q^{'}_{ij}, (0 \leq i \leq n, 0 \leq j \leq s) $ holds, then there is  
	$$
	\tilde{L}^{'} \tilde{P}_{\pi}=\tilde{P}_{\pi} Q^{'}.
	$$
	Furthermore $i m\left(\tilde{P}_{\pi}\right)$ is $\tilde{L}^{'}-$ invariant. 	
\end{lemma}

Similar to the proof of Theorem 1 and Corollary 1, we can get the following conclusion by Lemma 3, 
\begin{theorem}
	Suppose that $A_{i}=A_{j},B_{i}=B_{j}$ if nodes $i,j $ belong to the same cell and the matrix $Q^{'}$ exist, then every cell in the EP is trivial when system (3) is controllable.  
\end{theorem}
\begin{rmk}
	Theorem 3 extends the EP-based method to the controllability analysis of heterogeneous systems for the first time. The results show that the difference of node dynamic equations does not affect the relevant conclusions of EP if the nodes belong to different cells. From the 'equitable' defined by EP, two points can be called equitable points which not only means that the received signals are the same (the sum of the weight matrix is the same), but also their own dynamics should be the same. This deepens the understanding of EP.
\end{rmk}

\section{Union graph and its controllablity}
In this section, the controllability of the system is discussed from the perspective of the union graph, more specifically the relationship of controllability between union graph and system (2) is discussed, and a new criterion to identify the controllability of the union graph is obtained by using the concept of EP.
\begin{de}
The union graph of a collection of graphs $G_{m}=$ $\left\{V, E_{m}, \mathcal{A}_{m}\right\}, m=1, \ldots, t,$ is denoted by $\tilde{G}=\{V, \tilde{E}, \mathcal{A}\},$ where
$\tilde{E}=\bigcup_{m=1}^{t} E_{m},$ and $\tilde{\mathcal{A}}=\sum_{m=1}^{t} \mathcal{A}_{m}$.
\end{de}

Accordingly, union graphs can be described as
$$
\dot{x}(t)=\tilde{L}^{*} x(k)+\tilde{M} y(k)\label{1} \quad\quad\quad\quad (4)
$$
where
$$
\tilde{L}^{*}=\left(t\left[\begin{array}{cccc}
	A & 0 & \cdots & 0 \\
	0 & A & \cdots & 0 \\
	\vdots & \vdots & \ddots & \vdots \\
	0 & 0 & \cdots & A
\end{array}\right]-\left[\begin{array}{cccc}
	B K & 0 & \cdots & 0 \\
	0 & B K & \cdots & 0 \\
	\vdots & \vdots & \ddots & \vdots \\
	0 & 0 & \cdots & B K
\end{array}\right] L^{*}\right), \tilde{M}=\left(M\left[\begin{array}{cccc}
	C & 0 & \cdots & 0 \\
	0 & C & \cdots & 0 \\
	\vdots & \vdots & \ddots & \vdots\\
	0 & 0 & \cdots & C
\end{array}\right]\right)
$$
$L^{*}$ is the Laplacian matrix of the union graph.

In order to discuss the relationship between the controllability of system (2) and (4), the following lemma is introduced 
\begin{lemma}
System (2) is controllable if and only if the matrix
$$
\begin{array}{c}
	\tilde{W}_{c}=\left[\tilde{M}, \tilde{L}_{1} \tilde{M}, \ldots, \tilde{L}_{t} \tilde{M}, \tilde{L}_{1}^{2} \tilde{M}, \ldots, \tilde{L}_{1}^{dN-1} \tilde{M}\right., \left.\tilde{L}_{1}^{dN-1} \tilde{L}_{2} \tilde{M}, \ldots, \tilde{L}_{t}^{dN-1} \tilde{M}\right]
\end{array}
$$
is full row rank.
\end{lemma}
\begin{pf}
The proof is similar to Lemma 1 in \cite{8}, thus there is omitted.
\end{pf}
\begin{theorem}
	System (2) is controllable if the system (4) is controllable
\end{theorem}
\begin{pf}
From the Kalman rank criterion, Tthe following matrix is full row rank
$$
\begin{array}{c}
	\left[\tilde{M}, \tilde{L}^{*} \tilde{M}, \ldots, \tilde{L}^{* dN-1} \tilde{M}\right]
\end{array}
$$
which is expanded by using the properties of the union graph
$$
\begin{array}{c}
	\left[\tilde{M}, \left(\tilde{L}_{1}+\ldots+\tilde{L}_{t} \right)^{*} \tilde{M}, \ldots, \left(\tilde{L}_{1}+\ldots+\tilde{L}_{t} \right)^{* dN-1} \tilde{M}\right]
\end{array}
$$
Continue to expand the above matrix and add some column vectors, then
$$
\begin{array}{l}
	\left[\tilde{M}, \tilde{L}_{1} \tilde{M}+\cdots+\tilde{L}_{t} \tilde{M}, \tilde{L}_{2} \tilde{M}, \ldots, \tilde{L}_{t} \tilde{M}, \ldots, \tilde{L}_{1}^{d N-1} \tilde{M}\right. 
	\left.+\tilde{L}_{1}^{d N-2} \tilde{L}_{2} \tilde{M}+\cdots+\tilde{L}_{t}^{d N-1} \tilde{M}, \tilde{L}_{1}^{d N-1} \tilde{L}_{2} \tilde{M}, \ldots, \tilde{L}_{t}^{d N-1} \tilde{M}\right]
\end{array}
$$
Because the basic transformation of the matrix does not change the rank of the matrix, the above matrix can be transformed into
$$
\begin{array}{c}
	\left[\tilde{M}, \tilde{L}_{1} \tilde{M}, \ldots, \tilde{L}_{t} \tilde{M}, \tilde{L}_{1}^{2} \tilde{M}, \ldots, \tilde{L}_{1}^{dN-1} \tilde{M}\right., \left.\tilde{L}_{1}^{dN-1} \tilde{L}_{2} \tilde{M}, \ldots, \tilde{L}_{t}^{dN-1} \tilde{M}\right]
\end{array}
$$
which is full row rank. According to Lemma 4, system (2) is controllable.
\end{pf}
\begin{rmk}
Theorem 4 is only a sufficient condition, that is, when the system (4) is uncontrollable, the system (2) can also be controllable. See example 4 in Section 7 for more details.
\end{rmk}

From Theorem 4 and its examples, the controllability of system (2) cannot be determined when system (4) is uncontrollable only by above results. Therefore, we need to further analyze more specific cases. Inspired by Theorem 2, the following results are obtained
\begin{theorem}
Suppose the matrix $Q$ always exists for all $\tilde{L}_{i}$ corresponding to $G_{i}$, then the union graph is uncontrollable if any subgraph in the switching sequence contains nontrivial cells. 
\end{theorem}
\begin{pf}
If a proposition holds, so does its converse negative proposition. The converse negative proposition of Theorem 4 is that system (4) is uncontrollable if system (2) is uncontrollable. Theorem 2 shows that the dimension of the controllable subspace of system (2) is less than $dN$ if any subgraph in the switching sequence contains nontrivial cells. It means that the system (2) is uncontrollable, so we can get that the union graph is also uncontrollable. 
\end{pf}

\begin{rmk}
 Theorem 4 gives a special class of graphs. For this kind of graph, when system (4) is uncontrollable, system (2) is also uncontrollable. It avoids the uncertainty of previous results. At the same time, a new method to identify the controllability of the union graph is also given, which simplifies the process of identifying the controllability of the union graphs.
\end{rmk}

\section{Observability analysis}
Compared with controllability, observability is also a topic worth discussing, which is used to measure the ability of reconstructing the state of the whole network. Our previous discussion is based on the leader-follower structure. The leader is used to receive external input signals to control the whole network. Similarly, we can directly measure the state of some nodes, and then use these information to estimate the overall states of the network. Therefore, we call these nodes that can read the states directly as output points and the remaining nodes are called observed points. Accordingly, an architecture of output nodes-observed nodes is proposed to analyze the observability.

From the viewpoint of control, it is expected to find a controllable and observable system. For multi-agent systems, this requires that the controllability factor should be considered when selecting the output points. Therefore, in this paper, each leader in the network discussed above is selected as an output point, and the system can be described as
$$
\dot{x}(t)=\tilde{L} x(t)+\tilde{M} y(t)
\quad\dot{y}(t)=\tilde{M}^{T} x(t) \quad\quad\quad\quad (5)
$$
where
$$
\tilde{L}=\left[\begin{array}{cccc}
	A & 0 & \cdots & 0 \\
	0 & A & \cdots & 0 \\
	\vdots & \vdots & \ddots & \vdots \\
	0 & 0 & \cdots & A
\end{array}\right]-\left[\begin{array}{cccc}
	B K & 0 & \cdots & 0 \\
	0 & B K & \cdots & 0 \\
	\vdots & \vdots & \ddots & \vdots \\
	0 & 0 & \cdots & B K
\end{array}\right] L
$$
$$
\tilde{M}=M\left[\begin{array}{cccc}
	C & 0 & \cdots & 0 \\
	0 & C & \cdots & 0 \\
	\vdots & \vdots & \ddots & \vdots \\
	0 & 0 & \cdots & C
\end{array}\right].
$$

The observability of system (5) is equivalent to the controllability of its dual system. The dual system is represented as follows.
$$
\dot{x}(t)=\tilde{L}^{T} x(t)+\tilde{M} y(t),
\quad\dot{y}(t)=\tilde{M}^{T} x(t) \quad\quad\quad\quad (6)
$$
where $\tilde{L}^{T}$ is
$$
\begin{array}{c}
	\left[\begin{array}{cccc}
		A^{T} & 0 & \cdots & 0 \\
		0 & A^{T} & \cdots & 0 \\
		\vdots & \vdots & \ddots & \vdots \\
		0 & 0 & \cdots & A^{T}
	\end{array}\right]\\
	-L^{T}\left[\begin{array}{cccc}
		K^{T} B^{T} & 0 & \cdots & 0 \\
		0 & K^{T} B^{T} & \cdots & 0 \\
		\vdots & \vdots & \ddots & \vdots \\
		0 & 0 & \cdots & K^{T} B^{T}
	\end{array}\right]
\end{array}
$$
and $\tilde{M}$ is
$$
\tilde{M}=M\left[\begin{array}{cccc}
	C & 0 & \cdots & 0 \\
	0 & C & \cdots & 0 \\
	\vdots & \vdots & \ddots & \vdots \\
	0 & 0 & \cdots & C
\end{array}\right].
$$
And then we transform the observability problem of system (5) into the controllability problem of system (6). It is important to note that system (6) and system (1) are in different forms due to the existence of coefficient matrix, so the controllability of system (6) needs to be discussed.
\begin{lemma}
	$ \pi=\left\{V_{1}, V_{2}, \cdots, V_{k}\right\}$ is an EP for matrix-weighted graph $G$ and $P_{\pi}$ is the characteristic matrix.
	If the matrix
	$$
	Q^{*}=\left[\begin{array}{cccc}
		Q_{1*} & 0 & \cdots & 0 \\
		0 & Q_{1*} & \cdots & 0 \\
		\vdots & \vdots & \ddots & \vdots \\
		0 & 0 & \cdots & Q_{1*}
	\end{array}\right]-\left[\begin{array}{cccc}
		Q_{11*} & Q_{12*} & \cdots & Q_{1 s*} \\
		Q_{21*} & Q_{12*} & \cdots & Q_{2 s*} \\
		\vdots & \vdots & \ddots & \vdots \\
		Q_{n 1*} & Q_{n 2*} & \cdots & Q_{n s*}
	\end{array}\right]
	$$
	exists, which means that there are matrices $Q_{1*}, Q_{ij*}$ that make  $A^{T} C = C Q_{1*}$ and $L_{\pi ij} K^{T} B^{T} C = C Q_{ij*}, (0 \leq i \leq n, 0 \leq j \leq s)$ holds, then there is
	$$
	\tilde{L}^{T} \tilde{P}_{\pi}=\tilde{P}_{\pi} Q^{*}
	$$
	Furthermore $i m\left(\tilde{P}_{\pi}\right)$ is $\tilde{L}-$ invariant. 	
\end{lemma}
\begin{pf}
	The proof process is similar to Lemma 2, and thus is omitted here.
\end{pf}
\begin{rmk}
	Although Lemma 5 and Lemma 2 are similar in form, the restriction conditions of the coefficient matrix in the two cases are different, which is determined by the characteristics of the dual system.
\end{rmk}
\begin{theorem}
	Suppose that the matrix $Q^{*}$ exists, then system (5) is unobservable if the matrix weighted network contains nontrivial cells
\end{theorem}
\begin{pf}
	Similar to Theorem 1 and Corollary 1, we can get from Lemma 5 that the system (6) is uncontrollable if it contains the nontrivial cell. According to the properties of dual systems, system (5) is unobservable.
\end{pf}

If we take matrix $A$ as zero matrices, and matrices $B$,  $C$, and $K$ as identity matrices, then both matrices $Q$ and $Q^{*}$ exist according to Lemma 2 and Lemma 5. We call this case a first-order system. Choose and only choose all leaders as output points, then a result is obtained
\begin{co}
	System (5) is uncontrollable and unobservable if the network contains nontrivial cells under the first-order system.
\end{co}
\begin{rmk}
	Theorem 6 and Corollary 3 use the method of dual systems to discuss the influence of the existence of nontrivial cell on observability. Compared with the previous work of using EP to deal with observability of undirected and unweighted network \cite{43}, the situation considered in this paper is more extensive, and the method proposed in this paper is more general, which can be used to analyze the observability of multiagent systems in many cases.
\end{rmk}

\section{Examples}
In this section, we present several examples to demonstrate the relevant results in aboves.
\begin{figure}
	\centering
	\includegraphics[scale=.6]{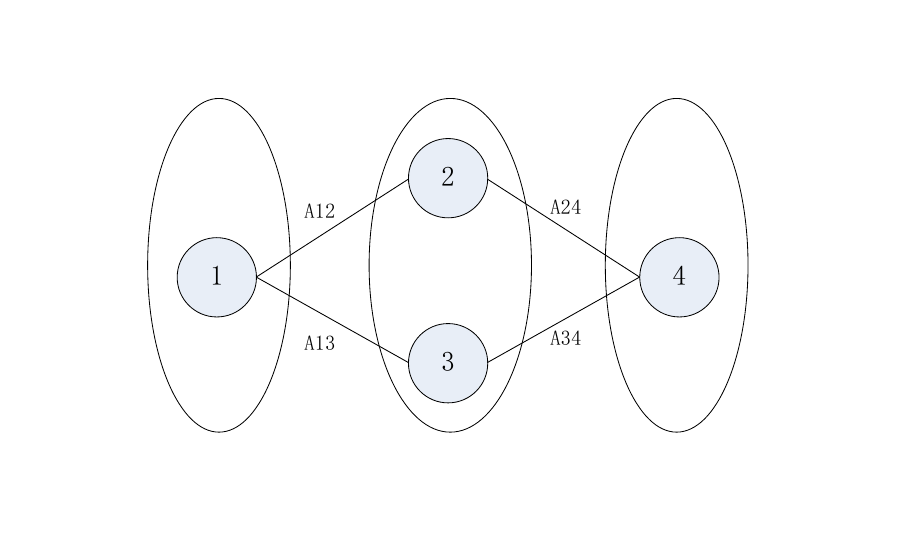}
	\caption{A matrix weighted signed graph with fixed topology including nontrival cell}
	\label{FIG:1}
\end{figure}
\begin{figure}
	\centering
	\includegraphics[scale=.6]{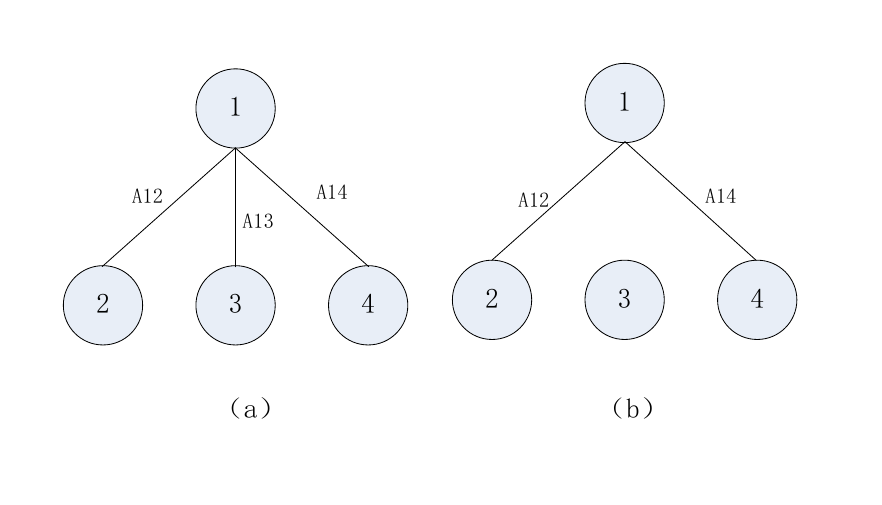}
	\caption{A matrix weighted signed graph with switching topology}
	\label{FIG:2}
\end{figure}
\begin{figure}
	\centering
	\includegraphics[scale=.6]{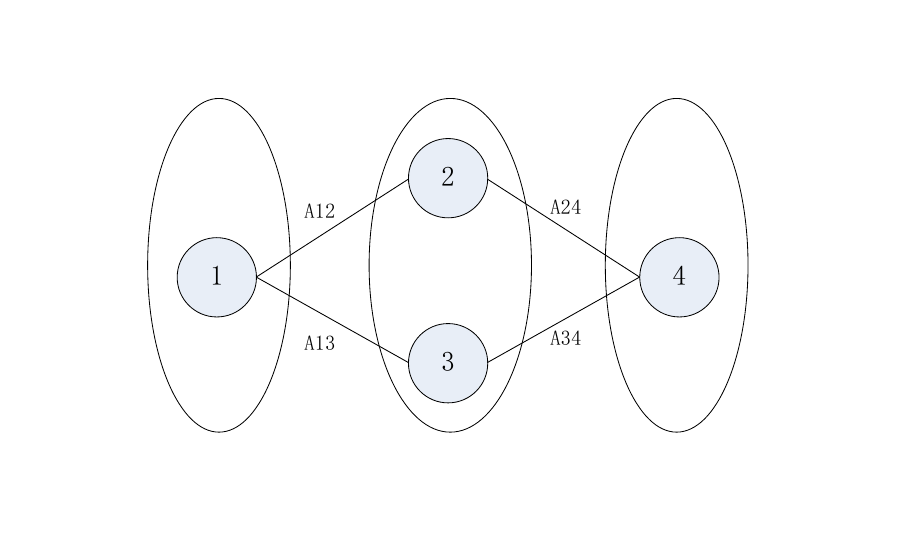}
	\caption{The case of heterogeneous system}
	\label{FIG:3}
\end{figure}
\begin{figure}
	\centering
	\includegraphics[scale=.6]{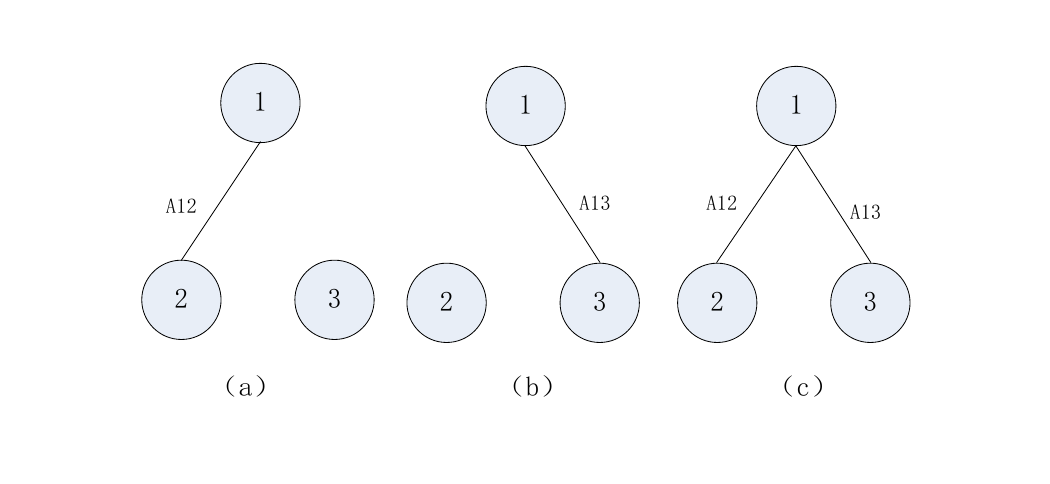}
	\caption{The case of union graph which is uncontrollable}
	\label{FIG:4}
\end{figure}	
\begin{figure}
	\centering
	\includegraphics[scale=.6]{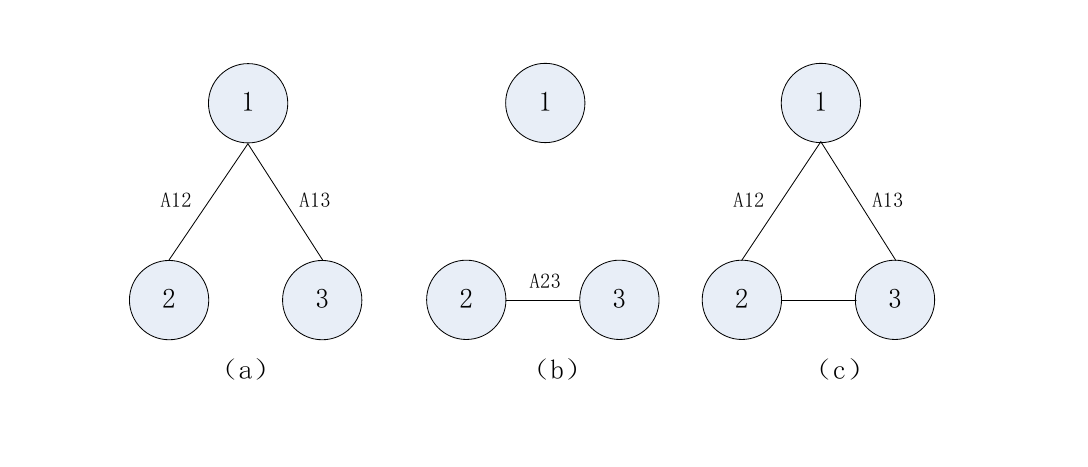}
	\caption{The case of subgraph including nontrival cell}
	\label{FIG:5}
\end{figure}
\begin{figure}
	\centering
	\includegraphics[scale=.6]{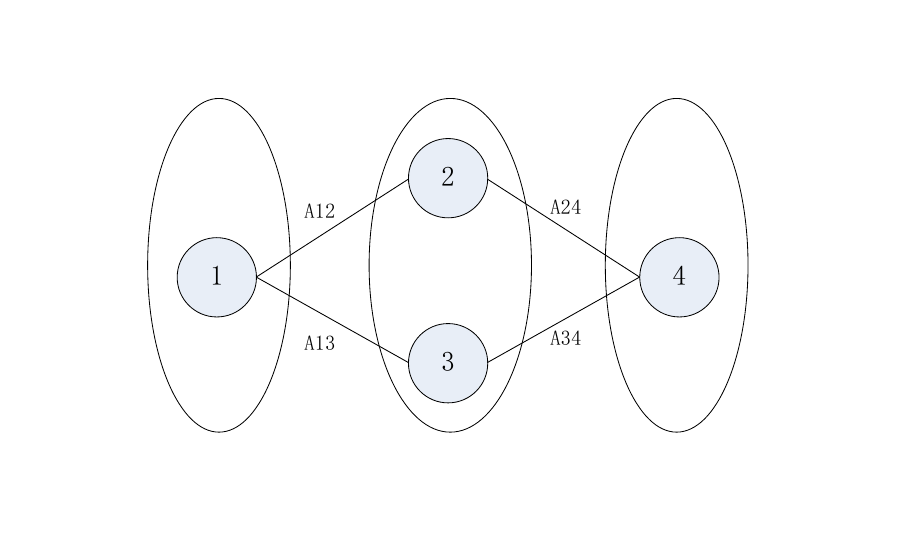}
	\caption{Observability of signed matrix weighted network}
	\label{FIG:6}
\end{figure}

\begin{exm}
	As showing in Figure 1, this is a network of four nodes. The coefficient matrix of the system is selected as
	$$
	A=K=\left(\begin{array}{ll}
		1 & 0 \\
		0 & 1
	\end{array}\right), B=C=\left(\begin{array}{ll}
		2 & 0 \\
		0 & 2
	\end{array}\right).
	$$
	The weights of the selected matrix are as follows
	$$
	\mathcal{A}_{12}=\mathcal{A}_{13}=\left(\begin{array}{ll}
		1 & 2 \\
		2 & 1
	\end{array}\right), \mathcal{A}_{24}=\mathcal{A}_{34}=\left(\begin{array}{ll}
		2 & 1 \\
		1 & 2
	\end{array}\right).
	$$
	Agent 1 is chosen as the signle leader. In view of the definition, the nodes can be divided into three cells $\pi=\left\{\left\{1\right\},\left\{2,3\right\},\left\{4\right\}\right\}$. The correcponding characteristic matrix is
	$$
	P=\left[\begin{array}{lll}
		1 & 0 & 0 \\
		0 & 1 & 0 \\
		0 & 1 & 0 \\
		0 & 0 & 1
	\end{array}\right]
	$$
	then we get
	$$
	L=\left(\begin{array}{cccc}
		\mathcal{A}_{12}+\mathcal{A}_{13} & -\mathcal{A}_{12} & -\mathcal{A}_{13} & 0 \\
		-\mathcal{A}_{12} & \mathcal{A}_{12}-\mathcal{A}_{24} & 0 & -\mathcal{A}_{24} \\
		-\mathcal{A}_{13} & 0 & \mathcal{A}_{13}-\mathcal{A}_{34} & -\mathcal{A}_{34} \\
		0 & -\mathcal{A}_{24} & -\mathcal{A}_{34} & -\mathcal{A}_{24}-\mathcal{A}_{34}
	\end{array}\right).
	$$
	After calculation, $rank(\tilde{L},\tilde{M})=6$. The system is uncontrollable. Theorem 1 is verified.
\end{exm}	
\begin{exm}
	The switching topology is shown in Figure 2. The coefficient matrix of the system is the same as the choice in Example 2, and the weights of the selected matrix are as follows
	$$
	\mathcal{A}_{12}=\mathcal{A}_{13}=\mathcal{A}_{14}=\left(\begin{array}{ll}
		1 & 2 \\
		2 & 1
	\end{array}\right).
	$$
	Agent 1 is chosen as the leader. Then, by definition, the nodes can be divided into two cells $\pi=\left\{\left\{1\right\},\left\{2,3,4\right\}\right\}$ for (a). For Fig. 2(b), the nodes can be divided into three cells $\pi=\left\{\left\{1\right\},\left\{2,3\right\},\left\{4\right\}\right\}$ for (b).
	$$
	L_{a}=\left[\begin{array}{cccc}
		\mathcal{A}_{12}+\mathcal{A}_{13}+\mathcal{A}_{14} & -\mathcal{A}_{12} & -\mathcal{A}_{13} & -\mathcal{A}_{14} \\
		-\mathcal{A}_{12} & \mathcal{A}_{12} & 0 & 0 \\
		-\mathcal{A}_{13} & 0 & \mathcal{A}_{13} & 0 \\
		-\mathcal{A}_{14} & 0 & 0 & \mathcal{A}_{14}
	\end{array}\right]
	$$
	$$
	L_{b}=\left[\begin{array}{cccc}
		\mathcal{A}_{12}+\mathcal{A}_{14} & -\mathcal{A}_{12} & 0 & -\mathcal{A}_{14} \\
		-\mathcal{A}_{12} & \mathcal{A}_{12} & 0 & 0 \\
		0 & 0 & 0 & 0 \\
		-\mathcal{A}_{14} & 0 & 0 & \mathcal{A}_{14}
	\end{array}\right].
	$$
	According to Theorem 2, the upper bound of the sequence of controllable subspaces is 4. After calculation, $\operatorname{dim}(\mathcal{C})=4$ too. Thus Theorem 2 is verified.
\end{exm}
\begin{exm}
	As showing in Figure 3, the coefficient matrix of the system is selected as
	$$
	A_{1}=\left(\begin{array}{ll}
		1 & 0 \\
		0 & 1
	\end{array}\right),	A_{2}=	A_{2}=\left(\begin{array}{ll}
		2 & 1 \\
		1 & 2
	\end{array}\right), A_{3}=\left(\begin{array}{ll}
		1 & 2 \\
		2 & 1
	\end{array}\right).
	$$
	the other coefficient matrix and weight matrix of the system are as Example 2. After calculation, $rank(\tilde{L}^{'},\tilde{M})=6$, The system (3) is uncontrollable. Theorem 3 is verified.	
\end{exm}
\begin{exm}
	The other switching topologies and their union graph are shown in Figure 4 where (c) is an union graph of (a) and (b). The coefficient matrix of the system is the same as the choice in Example 2, and the weight matrix is
	$$
	\mathcal{A}_{12}=\mathcal{A}_{13}=\left(\begin{array}{ll}
		1 & 2 \\
		2 & 1
	\end{array}\right).
	$$
	Agent 1 is chosen as the leader. For graphs in Fig. 3(a) and Fig. 3(b), there are
	$$
	L_{3 a}=\left[\begin{array}{ccc}
		\mathcal{A}_{12} & -\mathcal{A}_{12} & 0 \\
		-\mathcal{A}_{12} & \mathcal{A}_{12} & 0 \\
		0 & 0 & 0
	\end{array}\right]
	$$
	$$
	L_{3 b}=\left[\begin{array}{ccc}
		\mathcal{A}_{13} & 0 & -\mathcal{A}_{13} \\
		0 & 0 & 0 \\
		-\mathcal{A}_{13} & 0 & \mathcal{A}_{13}
	\end{array}\right].
	$$
	After calculation, $\operatorname{dim}(\mathcal{C})=6$. The system is controllable, but the union graph is uncontrollable. This example shows that the controllability of the system needs to be discussed separately when the union graph is uncontrollable.
\end{exm}
\begin{exm}
	The case of subgraph including nontrivial cell is shown in Figure 5, where (a), (b) are subgraphs and (c) is the union graph. The coefficient matrix of the system is selected as Example 2, and the weight matrices are as follows
	$$
	\mathcal{A}_{12}=\mathcal{A}_{13}=\mathcal{A}_{23}=\left(\begin{array}{ll}
		1 & 2 \\
		2 & 1
	\end{array}\right).
	$$
	Agent 1 is chosen as leader, and for Figure 3(a) and 3(b)
	$$
	L_{4 a}=\left[\begin{array}{ccc}
		\mathcal{A}_{12}+\mathcal{A}_{13} & -\mathcal{A}_{12} & -\mathcal{A}_{13} \\
		-\mathcal{A}_{12} & \mathcal{A}_{12} & 0 \\
		-\mathcal{A}_{13} & 0 &\mathcal{A}_{13}
	\end{array}\right]
	$$
	$$
	L_{4 b}=\left[\begin{array}{ccc}
		0 & 0 & 0 \\
		0 & \mathcal{A}_{23} & -\mathcal{A}_{23} \\
		0 & -\mathcal{A}_{23} & \mathcal{A}_{23}
	\end{array}\right].
	$$
	After calculation, $\operatorname{dim}(\mathcal{C})=4$, the system is uncontrollable, and the union graph is uncontrollable too. Theorem 5 is verified.
\end{exm}
\begin{exm}
	As showing in Figure 6, we choose node 1 as the output point and the coefficient matrix and weight matrix of the system are as Example 2. After calculation, $rank(\tilde{L}^{T},\tilde{M})=6$, The system (6) is uncontrollable, so system (5) is unobservable. Theorem 6 is verified.
\end{exm}

\section{Conclusion}
In this paper,  the controllability and observability of linear multi-agent systems over matrix-weighted signed networks were analyzed. Firstly, the definition of equitable partition of the matrix-weighted signed multi-agent system was given, and the upper bound of controllable subspace and a necessary condition of controllability were obtained by combining the restriction conditions of coefficient matrix and matrix weight for the case of fixed and switching topologies, respectively. The influence of different selection methods of coefficient matrices on the results was discussed. In particular, the equivalence of the controllability results based on equitable partition for first-order integrator systems and general linear systems under appropriate input matrix selection(it called complete control input) conditions was discussed. Secondly, for the case of heterogeneous systems, the upper bound of controllable subspace and the necessary condition of controllability were obtained when the dynamics of individuals in the same cell are the same. Thirdly, sufficient conditions for controllable and uncontrollable union graphs were obtained by taking advantage of the concept of switched systems and equitable partitions, respectively. Finally, a necessary condition of observability was obtained in terms of dual system and the constraints of the coefficient matrix, and the relationship between the observability and the controllability of the matrix-weighted signed multi-agent systems was discussed. For future work, we are interested in discussing the case that the weight matrix is asymmetric or other forms of weight (such as vector weight), to further improve the relevant results proposed in this paper.

\bibliographystyle{cas-model2-names}

\bibliography{ref}

\end{document}